\documentclass[10pt,a4paper]{article}

\usepackage[square,numbers]{natbib}
\usepackage{graphicx}
\usepackage{array,multirow}
\usepackage{makecell}
\usepackage{tablefootnote}
 \usepackage{amsmath}
\usepackage{amsfonts}
\usepackage{hyperref}
\usepackage[table,xcdraw]{xcolor}
\usepackage[margin=2.5cm]{geometry}

\usepackage[utf8]{inputenc}

\DeclareMathOperator*{\argmin}{arg\,min}

\begin{document}

\title{Structural engineering from an inverse problems perspective}

\author{
A. Gallet$^{1}$,  S. Rigby$^{1}$,  T. N. Tallman$^{2}$, X. Kong$^{3}$,\\  I. Hajirasouliha$^{1}$, A. Liew$^{1}$,  D. Liu$^{4}$, L. Chen$^{1}$,\\ A. Hauptmann$^{5,6}$, and D. Smyl$^{7}$ \\ 

\\

\small{$^{1}$Department of Civil and Structural Engineering, University of Sheffield}\\
\small{$^{2}$School of Aeronautics and Astronautics, Purdue University}\\
\small{$^{3}$Department of Physics and Engineering Science, Coastal Carolina University}\\
\small{$^{4}$School of Physical Sciences, University of Science and Technology of China}\\ 
\small{$^{5}$Research Unit of Mathematical Sciences, University of Oulu}\\
\small{$^{6}$Department of Computer Science, University College London}\\
\small{$^{6}$Department of Civil, Coastal, and Environmental Engineering, University of South Alabama}
}



\maketitle

\begin{abstract}
The field of structural engineering is vast, spanning areas from the design of new infrastructure to the assessment of existing infrastructure.
From the onset, traditional entry-level university courses teach students to analyse structural response given data including external forces, geometry, member sizes, restraint, etc. -- characterising a \emph{forward} problem (structural causalities $\to$ structural response).
Shortly thereafter, junior engineers are introduced to structural design where they aim to, for example, select an appropriate structural form for members based on design criteria, which is the \textit{inverse} of what they previously learned.
Similar inverse realisations also hold true in structural health monitoring and a number of structural engineering sub-fields (response $\to$ structural causalities).
In this light, we aim to demonstrate that many structural engineering sub-fields may be fundamentally or partially viewed as \textit{inverse problems} and thus benefit via the rich and established methodologies from the inverse problems community.
To this end, we conclude that the future of inverse problems in structural engineering is inexorably linked to engineering education and machine learning developments.
\end{abstract}



\maketitle

\pagebreak

\tableofcontents

\pagebreak

\section{Introduction} 
The estimation of structural response from loading and boundary conditions is a fundamental concept in structural analysis, from elementary Euler-Bernoulli beam theory to non-linear simulations involving complex structures subjected to extreme earthquake excitation.
In fact, numerical computation of structural response from known causalities characterises a \textit{forward} problem (cause $\to$ effect) and has rightly been the source of significant research since the advent of modern computing.
Amongst the myriad of computational frameworks, the finite element method (FEM) \cite{reddy2019introduction,surana2016finite,hughes2012finite,zienkiewicz1977finite}, finite difference \cite{virdi2006finite,liszka1980finite}, spectral element \cite{de2015comparison,kudela2007modelling}, and hybridisations \cite{e2017hybrid} have proven both widely applicable and successful over the years.
The implementation of such forward models have aided engineers in their ability to model, analyse, and design structures with arbitrary geometry and precision, contributing greatly to the presence of skyscrapers, supersonic aircraft, large cruise ships, and many more engineering examples.
In the near future, the pervasiveness of, for example, the FEM appears inevitable while its usefulness is unquestionable in structural engineering applications.

Pragmatically however, the final configuration of structural members is not known at the beginning of the design process; i.e. one iteration of a structural simulation is not generally sufficient in a real project.
This reality implies that the design of structures is an iterative process -- for example the identification of appropriately sized structural members, connections, and restraints (causalities) from design constraints, building codes, and environmental considerations (data).
As is often the case, the iterative design process is carried out initially using design tables, rules of thumb, handcrafted protocols, optimisation regimes, etc.
Nonetheless, this process is emphatically an \textit{inverse problem}, where an engineer is given data alongside design objectives and challenged to determine the appropriate structural configuration (causalities).

\begin{figure}[ht]
	\centering
	\includegraphics[width=9.0cm]{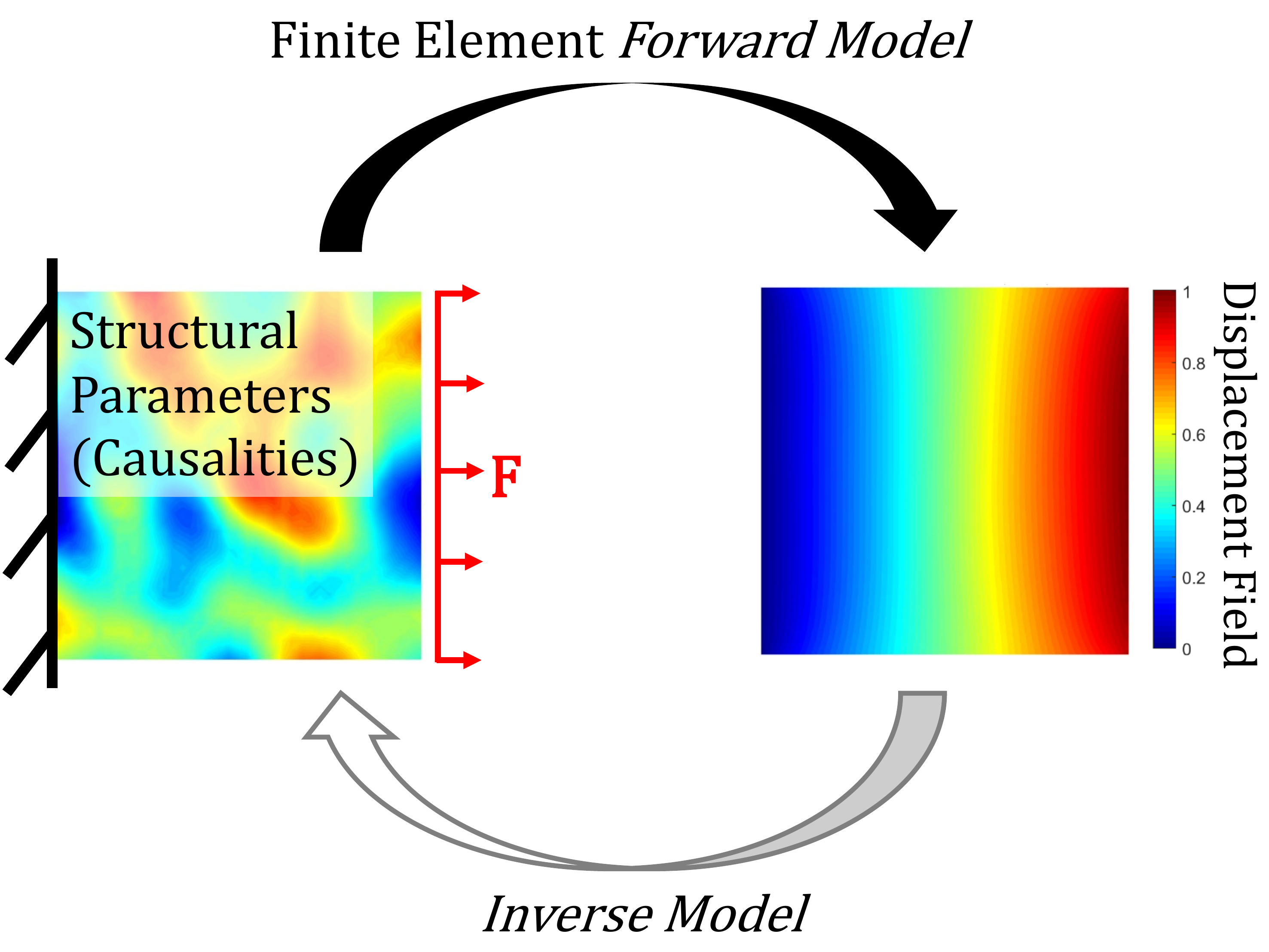}
	\caption{Schematic illustration depicting the forward and inverse problem relationship for a stretched elastic plate with randomised stiffness properties.  The forward finite element model inputs (causalities) are shown as the non-homogeneous stiffness properties while the model output is the displacement field. In contrast, the inverse model aims to estimate the stiffness properties \textit{given} the displacement field.} 
	\label{Strschem}
\end{figure}

Of course, the field of structural engineering is diverse, in which structural design is one of many sub-fields where inverse problems are applicable.
Perhaps a more straightforward implementation of inverse problems is structural health monitoring (SHM), where real-time (or near real-time) data is used in the prognosis of structural condition.
Indeed, the detection, localisation, prediction, and prognosis of potential damage processes is enabled in one of two ways, (a) via pattern recognition or (b) solving an inverse problem (or series of inverse problems) \cite{farrar2012structural}.
Moreover, certain non-destructive evaluation (NDE) approaches are also known to employ inverse problems (e.g. X-ray CT and emerging NDE approaches in academia) to asses the damage state of structural elements measured offline \cite{liu2003computational}.
For contextualisation, a schematic example contrasting the differences between forward and inverse models is provided in Fig. \ref{Strschem}; this well known problem is referred to as an \textit{inverse elasticity problem}.

These, and many more sub-fields of structural engineering {and structural mechanics \cite{turco2017tools}}, can not only be (fundamentally or partially) viewed as inverse problems, but, as we aim to illustrate herein, are benefited by the systematic approaches comprising the rich area of inverse problems.
Too often overlooked by structural engineers and structural researchers, the mathematics of inverse problems is an established field, ranging from classical statistical/Bayesian methodology \cite{kaipio2006} to cutting-edge implementation of deep neural networks \cite{adler2017solving}.
Moreover, while the present use of inverse problem methodologies in structural engineering is limited, its potential is immense across the expanse of the structural engineering sub-fields. 

In this paper, this potential will be discussed in detail and contextualised among a broad suite of existing inversion-based applications.
To begin, a clear description of inverse problems and methods will be detailed.
Following this, a review and discussion of inverse methodologies in modern structural engineering applications will be provided.
Inasmuch, the intent of this manuscript is to examine the following topics in structural engineering through the lens of inverse problems. 
We remark, however, that the forthcoming topical sections are \emph{not} intended to be exhaustive reviews, but rather, to provide substantiating evidence for the pervasiveness of inverse problems in structural engineering.
Lastly, realisations, overview, paths forward, and conclusions will be presented.


\section{Inverse problems, methods, and contemporary use}\label{invmath}
Traditionally, the field of inverse problems is concerned with the mathematical question of if and how one can determine the cause for certain measurements.  Despite being primarily mathematically oriented, the underlying questions always stem from relevant physics and engineering applications. This is especially true for one of the most prototypical inverse problems, the so-called Calder\'on problem \cite{calderon1980inverse} that asks: can one determine the conductivity of a body from electrical measurements at the boundary. In fact, this question arose during Calder\'on's time as a civil engineer, before he pursued an academic career in mathematics. In the following, we want to close the loop back to civil and structural engineering application that once motivated an entire field of mathematical studies, by utilising the insights gained in the last decades.

More generally speaking, inverse problems consist of finding the unknown characteristics of a structural system from some of the outputs, or measurements of that system. Most notably, this includes the above mentioned inverse conductivity problem in geophysics \cite{tarantola1982inverse} and engineering, but also includes a large field with applications in medical image reconstruction \cite{natterer2001mathematics,mueller2012}. Mathematically, such problems are ill-posed, broadly meaning that the parameters to be estimated $\theta$ are highly sensitive to changes in the measurement data $d$. The solution to the inverse problem involves estimating the parameter $\theta$ from a fixed set of measured data $d$, in contrast to the forward problem of computing $d$ from knowledge of the system parameter $\theta$. Specifically, this means given the forward model $U$, that models the system equations, we first formulate the underlying observation model
\begin{equation}\label{eqn:ForwProb}
    d=U(\theta) + \delta d.
\end{equation}
Here, $\delta d$ denotes an error term, modelling several sources of possible errors, such as inaccurate measurements or even inaccuracies in the model simulation.
The question that remains is: "how can we obtain $\theta$ from $d$ given the above relationship?"
Which we will call the reconstruction problem in the following. It is important to note that we can not simply invert the forward problem \eqref{eqn:ForwProb}, as the ill-posed nature implies that there can be either no or multiple solutions and additionally under inevitable measurement noise these solutions are not stable to compute by direct inversion. This ill-posedness of the inversion procedure constitutes the underlying paradigm of an \emph{inverse problem}.

In order to obtain stable reconstructions, we make use of a concept known as regularisation \cite{engl1996regularization}, which aims to assign a unique solution to each set of measurements in a stable manner, that means if the noise in the measurement vanishes, we would obtain the original system parameter. We can separate such stable reconstruction procedures into two primary classes: ones that compute a solution $\theta^*$ directly from measured data and those that iteratively aim to fit a solution by minimising a suitable cost functional. 
In the first case we aim to formulate an inverse mapping $U^\dagger$, such that 
\begin{equation}\label{eqn:directRecon}
    U^\dagger(d) \approx \theta.
\end{equation}
The primary problem in obtaining such direct inversion algorithms is that they can be highly dependent on the problem under consideration. Especially so, when the relationship between $d$ and $\theta$ is non-linear. Thus, obtaining such a mapping is a highly non-trivial task, but reveals much about the underlying problem characteristics and hence is a primary interest of mathematical research \cite{natterer2001mathematics,nachman1988reconstructions,nachman1996global}.

The second case, is a more principled approach that can be formulated for a large class of problems. The underlying premise is to reformulate the reconstruction problem as an optimisation problem. That is, we formulate a cost functional that measures how good our reconstruction fits the data while simultaneously enforcing some additional characteristics and acting as regularisation for the reconstruction process. Specifically, the reconstruction problem then writes as finding a minimiser of 
\begin{equation}\label{eqn:optiProb}
    \theta^* = \argmin_\theta \frac{1}{2}\|U(\theta) - d\|_2^2 + \alpha R(\theta).
\end{equation}
Here, the first term enforces that reconstructions fit the data, whereas the second term is the so-called regularisation term. As discussed previously, this regularisation term is necessary when dealing with inverse problems, as it prevents a solution from over-fitting the measurement noise.
Importantly, by incorporating prior knowledge in the design of $R$ \cite{tikhonov1977,RudinOsherFatemi1992}, we effectively choose preferred solutions and overcome the problem of non-uniqueness. Finally, the parameter $\alpha>0$ balances both terms and depends on the noise amplitude. Solutions to \eqref{eqn:optiProb} are computed by suitable optimisation schemes, for which repeated evaluation of the forward model $U$ will be necessary. Consequently, computing solutions to \eqref{eqn:optiProb} can be computationally highly expensive, if the evaluation of $U$ is expensive. Thus, for nonlinear problems fast converging algorithms, such as Gauss-Newton or related methods \cite{haber2000optimization}, are preferred.

Lastly, with the recent rising popularity of data-driven methods, researchers have designed computationally more efficient ways to address the reconstruction problem \cite{arridge2019solving}. Such data-driven approaches are often inspired by classical reconstruction algorithms discussed above. For instance, one can replace the direct reconstruction operator or parts of it, with a data-driven component, typically consisting of a neural network \cite{jin2017deep,kang2006,hamilton2018deep}. Given a set of informative reference data, one can then learn a suitable mapping mimicking \eqref{eqn:directRecon}. Alternatively, researchers have investigated the possibility to improve the iterative process to compute solutions to \eqref{eqn:optiProb}, by replacing parts in the optimisation algorithm with learned components \cite{adler2017solving,adler2018learned,hammernik2018learning,hauptmann2018model,smyl2021efficient}, or entirely building network architectures motivated by the iterative solution process \cite{zhang2018ista,monga2021algorithm}. {The use of data-driven methods is expanded upon in \S\ref{MLInversion}.}

\section{Structural design as an inverse problem}\label{invdes}

\subsection{Demarcating structural analysis and design}
Prior to World War II, engineering higher education was originally focused on the art and practice of engineering design \cite{crawley2007rethinking}. By the 1960s however, due to the success of science-based ventures such as the Manhattan Project and the rise of government-sponsored research grants that severed the link between academia and industry, engineering science became the main field of research and teaching at universities \cite{Seely1999}. This research led to the development of powerful structural analysis techniques, yet also left engineering graduates with a noticeable loss of practical engineering skills \cite{Mann2016}. In the 1990s and early 2000s, there was a push to introduce capstone design projects in university engineering degrees to address this shortcoming, the success of which is disputed \cite{Coleman2020}.

There is a broad agreement within the literature that analysis and design are two distinct activities. Structural analysis, which falls into the field of engineering science, is primarily concerned with establishing \textit{knowledge-that} explains the world, whereas design is concerned with \textit{knowledge-how} something works \cite{Bulleit2015}. Design is often characterised as being ill-structured \cite{Simon1973}, open-ended \cite{Dym2006} or even “wicked” \cite{Rittel1973}, qualities which do not necessarily lend it to science-based research and helps motivate the search for a different paradigm which more adequately addresses its true nature.

One perspective to account for these differences is to recognise structural analysis and structural design as being two different types of problems. Although structural analysis is typically seen as a sub-field of design to validate and justify the adequacy of structural elements \cite{Mason2018}, we like to advocate the view that structural design is in fact an \textit{inverse-problem}, with structural analysis forming the related \textit{forward-problem}.

\begin{figure}[h!]
	\centering
	\includegraphics[width=13cm]{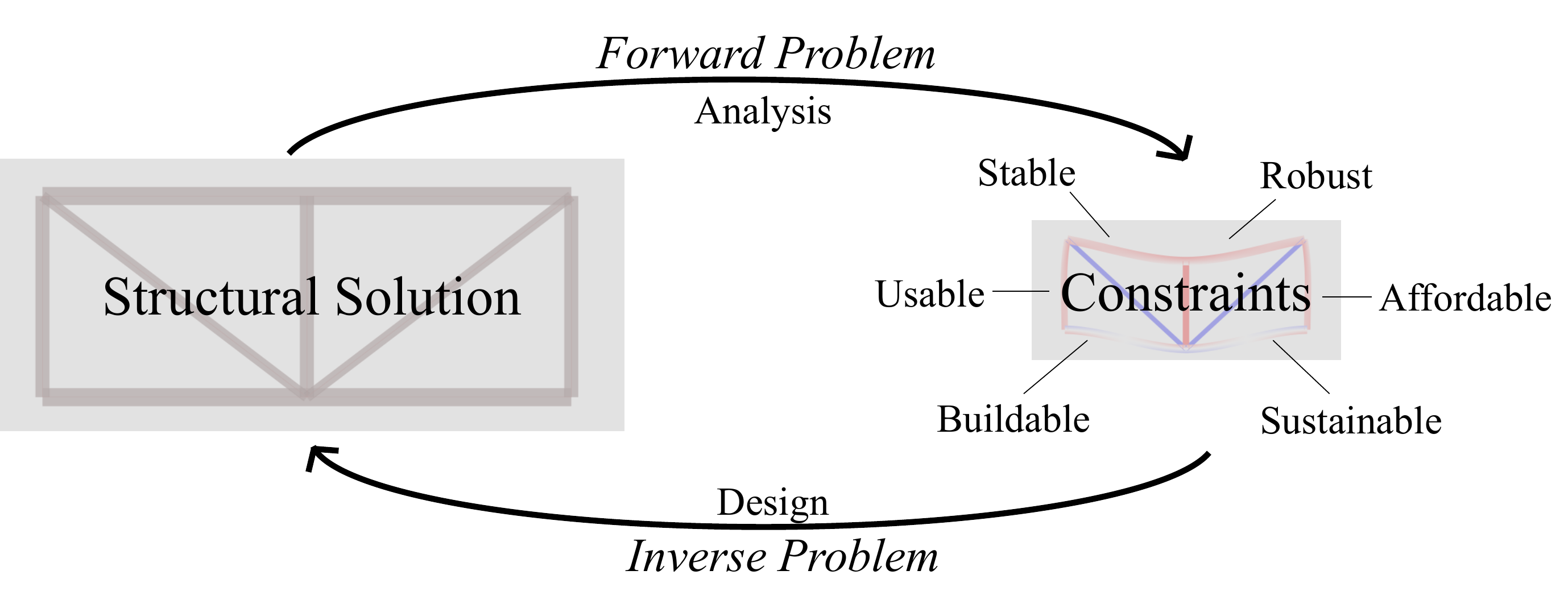}
	\caption{Relationship between structural analysis and design, in which design is the inverse problem of evaluating a suitable structural solution \textit{given} a set of constraints. Some constraints (such as stability constraints) render themselves more easily to quantitative treatments.}
	\label{fig:DesInv}
\end{figure}

Unlike typical inverse problems such as the one shown in Fig. \ref{Strschem}, where physically measured data (the displacement field) is used to identify the causalities (stiffness properties), in structural design we are dealing with a theoretical construct, or “theoretically” measured data. This data contains the set of complex design constraints which need to be adhered to, such as ultimate (ULS), serviceability limit states (SLS), sustainability and constructability. From this perspective, structural design can be seen as the process of arriving from a specific set of constraints to a viable structural solution, with analysis being the process of checking if the proposed structural solution adheres to those constraints as shown in Fig. \ref{fig:DesInv}.

This perspective is evidenced by realising that the key features of inverse problems, specifically their ill-posed nature as explained in \S\ref{invmath}, manifest themselves in design. Their shared characteristics include being unstable, non-unique and sometimes un-solvable problems \cite{Alifanov1983}. Notice how these qualities were alluded to as “ill-structured”, “open-ended” and “wicked” by previous researchers, yet to the best of our knowledge, this is the first time that literature has attempted to view structural design as an inverse problem.

{We suggest that the application of the inverse problem perspective in design gives rise to the idea of a ``design model''. Similar to how an ``analysis model'' allows us to evaluate the action effects of a given structural solution, solving the forward problem shown in Fig. \ref{fig:DesInv}, a ``design model'' would generate structural solutions which adhere to the given set of user defined constraints, solving the inverse problem. This perspective sheds light on the characteristics which design problems directly share with other inverse problems, such as the non-uniqueness property. Multiple viable solutions often exist to any given design brief, one example of which are the various truss-archetypes for bridge designs to span a similar distance. In order to create a viable design model, it would be necessary to provide some form of regularisation, examples of which are provided in \S\ref{shm} in the context of structural health monitoring, which effectively encourage design models to choose preferred solutions based on prior knowledge. This could be achieved by constraining the solution space to a sub-set of structural systems, cross-sections or materials based on the specific constraints provided. Other properties of inverse problems which arise in design are shown in Table \ref{table:DesInvQual} and a comprehensive understanding of the application of the inverse problem perspective warrants further research using specific design examples.}

\begin{table}[ht!]
    \centering
    \scriptsize
    \begin{tabular}{>{\centering\arraybackslash}m{0.15\textwidth} !{\vrule width 1pt} >{\centering\arraybackslash}m{0.35\textwidth} | >{\centering\arraybackslash}m{0.35\textwidth}}
        \textbf{Ill-posed inverse problem characteristics} & \textbf{General description} & \textbf{Examples in structural design} \\
        \Xhline{1pt}
        Unstable\tablefootnote{In reference to the solution space, not structural instability (buckling), see \cite{Alifanov1983}.} & Small changes in constraints can lead to large changes in the solution space & Impact of grid-spacing on the selection of appropriate floor options \\
        \hline
        Non-unique & For a given set of constraints, multiple solutions exist & Various truss archetypes which can span the same distance \\
        \hline
        Unsolvable & Lack of realistic constraining data impedes the search for an adequate solution & Current lack of appropriate structural materials for space elevator designs \\
        \Xhline{1pt}
    \end{tabular}
    \caption{Overview of ill-posed inverse problem features in the context of structural design.}
	\label{table:DesInvQual}
\end{table}

We note that the design problem can be also formulated as an optimal control problem \cite{yang1975application}, where the optimal design parameters are thought to be found as minimiser of a penalty function while satisfying the system equations. Whereas, the optimal control approach can provide an effective way to solve complicated design problems, it falls short in accounting for uncertainties or inaccuracies in the forward model, and especially the link between measurement and system parameters. We believe that the here lies the strength of the inverse problems viewpoint, that offers a rich interpretation and link between system parameters and measurements given as the forward model. This is not only promising for the optimisation task in the design process, but also for new ways to approach the modelling of the forward problem.

\subsection{The link to structural optimisation}
As explained in \S\ref{invmath}, inverse problems can effectively be solved iteratively \cite{Wise2016}. Effectively, this process involves making an estimate of the structural solution based on experience and design heuristics (such as simple rules of thumbs), checking those estimates through an analysis (forward) model and updating the model if necessary; in other words, this approach is characterised by creating an optimisation problem. Some examples of forward-driven optimisation models used in structural design include: optimising the deformed shape of flexible formwork structures to predefined target geometries \cite{ROMBOUTS2021112352,LIEW2018128}, best-fit geometry optimisation of thrust networks in the design of shell structures \cite{VanMele2014}, and finding the optimal structural forms for long-span bridges \cite{fairclough2018theoretically}, gridshells \cite{LIEW20201845}, trusses \cite{HEGILBERT}, portal frames \cite{mckinstray2015optimal} and structural sections \cite{YE201855}. 

A key theme in these research works is that the structural geometry or member proportions are not initially defined, but rather are form-found or discovered in the process, based on the defined loading, boundary conditions and objectives. Often these discovered structural forms may lead to step change benefits in terms of performance or reduced material usage, as they are unbiased by our preconceived perception of what a "good" structure is, or by what we currently design and build in the construction industry with standard template solutions. It is also true for many cases, that a solution may not even exist, forcing us to accept a closest best fit solution, or it can be that an optimal solution may have multiple candidates by virtue of the structure's static indeterminacy. A common problem with a forward (sometimes brute-force) optimisation approach, can be lengthy computational times for structural analyses in the objective functions, and the sheer size of the design (and hence optimisation) search space, stemming from many wide ranging input design variables. While fast and globally convergent convex optimisation programs can be set-up, many practical structural engineering design problems are inherently non-linear in nature, forcing a slower approach that uses local searches or heuristic methods with no guarantee of a global optimum. This is a current challenge faced in solving inverse problems iteratively with forward models.

\subsection{Implications of treating structural designs as inverse problems}
This inverse problem perspective has various implications. Firstly, it should emphasise that design problems ask a different question than those related to analysis. An analysis model solves the forward-problem, and answers questions such as ``what is the ULS utilisation ratio of a particular beam system for a specific load, with the following specific cross-sections and support conditions?''. An appropriate design model would ask the reverse: ``what is the group of cross-sections and support-conditions which ensures a utilisation ratio of less than 1.0 for this particular beam system to carry this particular load?''. In this formulation, the magnitude of loading and the utilisation ratio serve as the design constraints (both are known ``data''). Whilst engineering science has produced sophisticated analysis models, the research of such ``design models'' is lacking in academia.

Secondly, the inverse problem perspective sheds insight on the possibility of using data-driven approaches as opposed to more typical optimisation techniques. The rise of machine learning and deep neural networks could be used for the development of such ``design models'', which focus on directly identifying a set of structural solutions from a given set of constraints using learned data. Such models could address some of the challenges faced in design due to tight deadlines that force early design decisions, whose full implications might only be realised at the detailed design stage where changes become cost and time prohibitive. This can lead to structural solutions that are difficult to build, have poor sustainability metrics, and be costly to engineer and fabricate. If one instead considers the design process in an inverse manner, rooting firmly first at the end goal, it could be possible to reduce project risk and pick more effective structures by appreciating many solutions to the brief from the onset.
{This was already alluded to in \S\ref{invmath} and is expanded upon in section \S\ref{MLInversion}}.

Lastly, the ultimate benefit of an inverse problem perspective is that it helps to clearly distinguish between analysis and design procedures and provides academia with an adequate framework to contextualise \textit{design model} research. Engineering academia has been dominated by the engineering science perspective \cite{Koen2003}, predominantly choosing to research and teach forward problems. One of the uncomfortable implications of this view is that engineers from over 150 years ago, which trained primarily in the art and practice of engineering design, may in fact have been better “inverse-problem solvers” than academics and graduating engineers of today (who are stronger in solving the well-structured forward problems) \cite{Addis1990}. This might help swing the pendulum away from focusing exclusively on forward models (analysis) towards a more stable equilibrium with inverse models (design) by acknowledging the existence of these two related, albeit distinctively different types of problems. 
{The use of inverse problem and inverse analysis in a related field to structural design, notably blast engineering, is discussed next.}


\section{Extreme loads on structures}\label{extreme_loads}
\subsection{Blast loading and inverse analysis}

High-rate dynamic loads can arise from events such as earthquakes, wind, tidal waves, impacts, and accidental or malicious explosions. Here, the imparted load may be comparable or several times larger than the strength of the material it is acting on, it can be applied and removed in sub-second durations, and is often highly localised. Accordingly, the notion of static structural design according to a pre-determined distribution of stresses and strains may not be appropriate, and instead the designer must consider energy balance, non-linear analyses, and deformation modes for which there is no equivalent static counterpart.

Blast loading is undoubtedly one of the most aggressive forms of dynamic loading. When an explosive detonates it undergoes a violent and self-perpetuating exothermic chemical reaction, releasing energy through the breaking of inter-molecular bonds during oxidation \cite{Cormie2009}. The explosive material is converted into a high pressure (10--30~GPa), high temperature (3000--4000$^\circ$C) gas which violently expands, displacing the surrounding air at supersonic velocities (6--8~km/s). This displacement causes a shock wave to form in the air, termed a blast wave, which eventually detaches from the expanding detonation product `fireball' and continues to propagate into the surrounding air, decreasing in pressure and density as it expands. 

When a blast wave encounters an obstacle some distance from the source it will impart momentum as the air is momentarily (either fully or partially, depending on the compliance of the obstacle) brought to rest at the air/obstacle interface. Prediction of blast even in the most simple settings is a considerable challenge to the scientific community. This becomes an increasingly complex and multi-faceted problem when considering issues such as: obstacle orientation; proximity of the obstacle to the source and additional momentum transfer from fireball impingement; secondary combustion effects either at the air/obstacle interface or in late-time due to partial or full confinement of the explosive products; the presence of mitigating or blast-enhancing materials (soil, reactive munitions, etc.).

Real-world blast events are highly uncertain, and the need for inverse analysis is clear: it is very rare that the exact size, shape, composition, and location of an explosive device is known \emph{a priori}. Instead, information relating to the \emph{cause} of an explosion should be estimated, within reason, from the more readily observed \emph{effects}, i.e.\ the magnitude and severity of structural damage to surrounding buildings, and cratering of the ground surface. Whilst inverse analysis is well-established for practical post-event assessment of explosions\textemdash and has been used to determine the size/location of blast events through forensic investigation of social media videos \cite{Rigby2020Beirut} or numerical modelling correlating structural damage \cite{Ambrosini2005,VanDerVoort2015}\textemdash the use of inverse modelling in an academic context is yet to be exploited fully. In the former, order-of-magnitude estimates are typically deemed sufficient, whereas the latter requires repeatable, precise measurements and high levels of experimental control.

The lack of robust yet high-fidelity experimental techniques has stifled academic research into close-in blast for some time. Close-in blast is typically defined as the region within approximately 20 radii from the charge centre, where loading in this region is characterised by a near-instantaneous rise to peak pressure in the order of 100--1000~MPa, followed by a rapid decay to ambient conditions typically occurring in sub-ms. Subsequent structural response may reach a peak value in the order of 10-100~mm and, whilst this may occur orders of magnitude slower than load application, deformation cycles are still generally within ms durations.

Recently, two notable advancements have been made in experimental characterisation of close-in blast and structural response. In the first of these, researchers at the University of Sheffield (UoS), UK, developed a large scale apparatus for the spatial and temporal measurement of blast pressures from close-in explosions \cite{Clarke2015,Rigby2015near}. In the second, researchers at the University of Cape Town (UCT), South Africa, adapted the well-known digital image correlation (DIC) technique to measure the transient response of the rear-face of blast loaded plates \cite{Curry2017}. These two techniques were combined in a recent study \cite{Rigby2019ExpMech,Rigby2019IJIE} where, for the first time, detailed loading maps \emph{and} temporal structural response profiles were developed independently, in a single-blind study, for identical (scaled) experimental set ups.

\subsection{Proof-of-concept experimental studies}

Of the experiments performed in \cite{Rigby2019ExpMech}, 12 are relevant to the notion of inverse analysis of blast loading and structural response, and will be discussed here. Six tests were performed with spherical explosive charges; three at UoS measuring blast loading and three at UCT measuring structural response. In the UoS spherical tests, 100~g PE4 charges were located at 55.4~mm clear distance from the centre of a nominally rigid target, on which the reflected blast pressures were measured. In the UCT spherical charge tests, 50~g PE4 charges were located at 44.0~mm clear stand-off distance from the centre of the flexible target plates: 300~mm diameter, 3~mm thick, Domex 355MC steel plates, fully clamped around the periphery. The plate response was filmed using a pair of stereo high speed video cameras and DIC was used to determine transient plate response. The two test series can be expressed at the same scale using well-known geometric/cube-root scaling laws. Here, it is assumed that the flexible targets deform on timescales orders of magnitude longer than loading application, and therefore differences between the loaded imparted to the rigid and flexible plates are negligible.

Specific impulse is given as the integral of pressure with respect to time. Numerical integration of the UoS pressure histories (at various distances from the centre of the plate) yields \emph{directly-measured} specific impulse distributions. The first few frames of the UCT tests were used to determine initial velocity distributions of the plate, from which imparted specific impulse could be \emph{inferred} through localised conservation of momentum: $i(x) = v(x)\rho t$, where $i$ is specific impulse, $x$ is distance from the plate centre, $v$ is out-of-plane velocity of the plate, $\rho$ is density (7830~kg/m$^3$ for Domex 355MC), and $t$ is plate thickness.

\begin{figure}[h!]
	\centering
	\includegraphics[width=\textwidth]{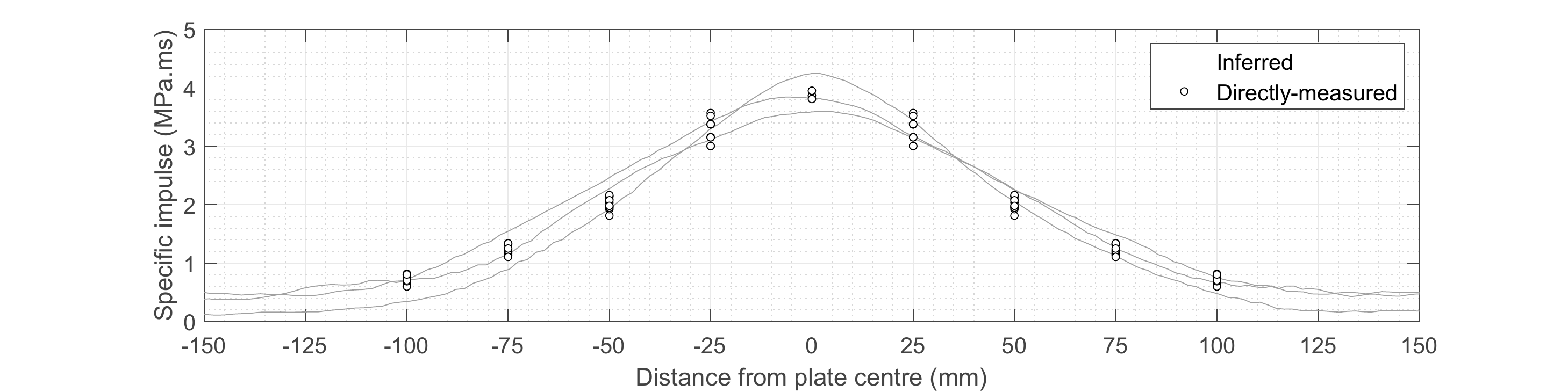}
	\caption{Directly-measured (UoS) and inferred (UCT) specific impulse distributions from studies of blast loading and plate deformation following detonation of spherical explosives, expressed at 100~g (UoS) scale, after \cite{Rigby2019ExpMech}}
	\label{fig:BlastSPH}
\end{figure}

The results for the spherical tests are shown in Fig.~\ref{fig:BlastSPH}. The full-field inferred specific impulse distributions are in close agreement with the discreet, directly measured values, and both measurements form a tight banding in an approximate Gaussian distribution \cite{Pannell2021}. Not only does this indicate a high level of test-to-test repeatability for each method, but demonstrates that the two methods are measuring the same underlying phenomena, albeit in entirely different ways. Thus, it can be said that an imparted impulse will result in an initial velocity uptake which is directly proportional, and therefore measurement of one allows for the other to be determined. This proves the concept of using plate deformation under blast loads in an inverse approach\textemdash namely that from knowledge of plate deformation one may be able to determine the imparted load\textemdash and provides physical verification of the inverse approach developed by \cite{Xu2010,Xu2011}.

In addition to the spherical charge tests, six tests were performed in \cite{Rigby2019ExpMech} using squat cylinders (height:diameter of 1:3). Such charges are known to produce a more concentrated load, with the fireball propagating at higher velocities along the axis of the charge \cite{Langran-Wheeler2021}. This accelerates the growth and emergence of surface instabilities \cite{Tyas2016}, which gives rise to a more variable loading distribution \cite{Rigby2020}. A key research question in this study was: "will a more variable loading result in more variable structural response?" In the UoS cylindrical tests, 78~g PE4 charges were located at 168.0~mm clear distance from the centre of the target, and in the UCT cylindrical tests, 50~g PE4 charges were located at 145.0~mm clear distance from the centre of the target. Again, specific impulse distributions were both directly-measured and inferred from plate response respectively, and the experiments were expressed at the same scale using common scaling laws.

\begin{figure}[ht!]
	\centering
	\includegraphics[width=\textwidth]{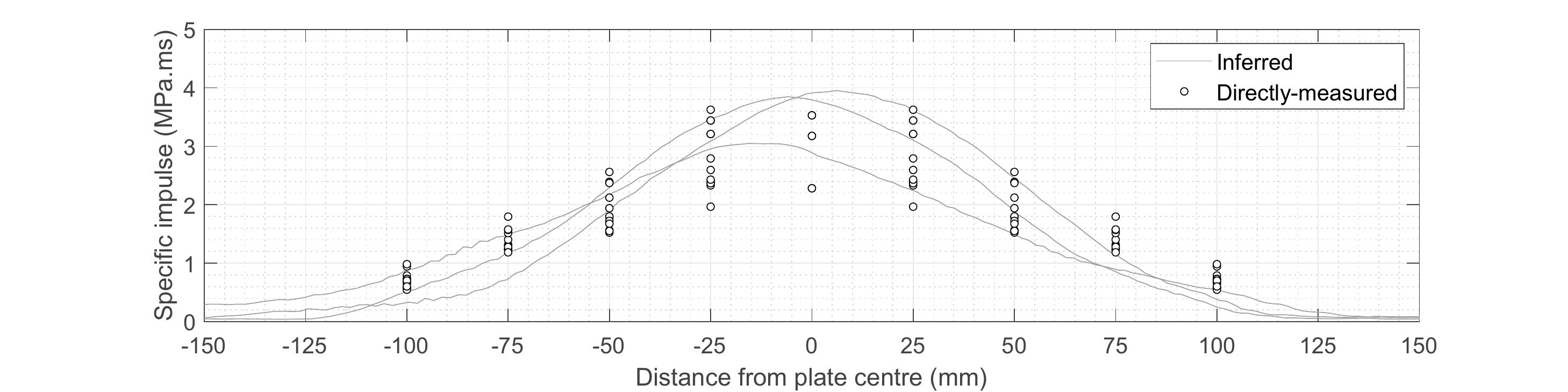}
	\caption{Directly-measured (UoS) and inferred (UCT) specific impulse distributions from studies of blast loading and plate deformation following detonation of cylindrical explosives, expressed at 78~g (UoS) scale, after \cite{Rigby2019ExpMech}}
	\label{fig:BlastCYL}
\end{figure}

The results for the cylindrical tests are shown in Fig.~\ref{fig:BlastCYL}. Whilst the two methods again show good agreement, the results can be seen to form a much wider spread. In contrast to the spherical tests, where peak specific impulse was seen to consistently act in the plate centre for all tests, here the peak value is often up to 25~mm from the plate centre (approximately equal to the charge radius), in both the directly-measured and inferred values. The inferred values are generally bounded by the directly-measured values, which suggests that this spread is indeed a genuine feature caused by application of more variable load.

\subsection{Outlook}
The aforementioned studies \cite{Rigby2019ExpMech,Rigby2019IJIE} have provided a firm physical basis for inverse analysis in the context of extreme loading and structural response. The results have clear implications for the future of research in this area. Namely, it has been demonstrated that not only can inverse analysis provide excellent predictions of blast loading in repeatable, well-controlled situations (as with the spherical tests in \cite{Rigby2019ExpMech}), but structural response measurements are potentially sensitive enough to detect localised variations in loading (as with the cylindrical tests in \cite{Rigby2019ExpMech}). This is particular important in situations where a highly variable loading might be expected (e.g.\ from explosives buried in well-graded soils \cite{Rigby2018CoBL}), but statistical variations cannot be determined in a robust sense when using direct measurements (note the discrete nature of the direct measurements in this study, compared to the effectively continuous nature of the inferred measurements). This technique may permit, through inverse methods, fundamental scientific studies of complex mechanisms governing blast loading following close-in detonation of explosive charges in situations where previous research has not yet been possible.

\section{Structural health monitoring}\label{shm}
\subsection{Background on inverse methods in SHM}

{Unlike the fields of structural design or blast engineering, it is well known that} inverse problems are deeply connected to SHM.
In fact, computerised damage detection can generally be recognised as either a pattern recognition or inverse problem \cite{farrar2012structural,friswell2007damage,kerschen2006past} where unknown or uncertain parameters (causalities) are estimated via quasi-static or dynamic structural response data.
Among the numerous inversion-based approaches, FEM updating methodologies are among the most pervasive \cite{sinha2002simplified,sinou2009review,farhat1993updating,ricles1992damage}.
Much of the noted popularity is owed to the flexibility of the FEM for comparing modal parameters to an undamaged state and in compensating model errors.

Meanwhile, analytical \cite{lin2004direct}, wavelet \cite{hester2012wavelet}, fractal \cite{hadjileontiadis2005fractal}, fuzzy system \cite{pawar2003genetic},  Kalman filter \cite{lourens2012augmented}, chaotic interrogation \cite{nichols2003structural}, shape function \cite{liu2014novel}, and particle filtering \cite{nasrellah2010particle} approaches, among many others, have proven successful in uncertain inverse parameter identification applications within SHM.
Broadly speaking, inverse SHM approaches can be grouped as either deterministic or probabilistic \cite{ebrahimian2018bayesian} -- which is also generally the case in classical inverse problems \cite{kaipio2006}.
In the latter case, estimation of uncertain SHM parameters takes the form of a probabilistic term, for example a value with an associated certainty, a probability itself, etc.

Irrespective of the computational approach used in damage detection, two key realisations affect the efficacy of inversion methodologies: (a) a baseline is generally needed to detect/quantify damage \cite{Farrar2007} and (b) the presence of damage inherently influences the linearity of structural behaviour \cite{worden2008review}.
In addressing (a), reference-free or baseline-free frameworks have been introduced \cite{huang2017baseline,sohn2007combination,santos2015baseline,anton2009reference} via the introduction of either assumptions on the reference state, implementation of prior physics knowledge, or probabilistic regimes.
On the other hand, non-linearity in the structural response can either act as a corrupting entity when linear forward models are used (e.g. unacceptable forward model error) or used as an advantage when properly leveraged.
Regarding the latter, as noted in \cite{worden2008review}, methods based on non-linear indicators, dynamical systems theories, and non-linear systems identification approaches can be used to aid or enrich the damage identification process; such a conclusion can also be extended to the pure usage of inverse approaches in damage detection. 

In the past thirty years, implementation of inversion-based damage detection methods in SHM has steadily increased.
This is due to both the increase in inverse problems know-how and computational resources.
Yet, since the emergence of contemporary machine learning, the ability to solve problems deemed previously intractable has exponentially increased opportunities in this area.
For example, in many cases, forward models may not be available or are too computationally expensive, sufficient non-linearity may exist to effectively model the desired physics,  errors in highly reduced models may be excessive, the ability to compute model gradients may be overly expensive, etc.
Moreover, the ability of classification networks to readily classify important variables such as the probability and/or severity of damage from structural data is intuitively appealing and pragmatically useful.
In the following, we will provide contemporary examples highlighting the use of both classical and machine learning based SHM inverse approaches.

\subsection{Static inverse problems in SHM}\label{static_inverse_problems_in_shm}
Incorporation of discrete static (or quasi-static) data measured from structures into SHM frameworks is well established.
For example, a number of sources, including corrosion, relative humidity, fibre-optic, topography, laser, potentiometer, strain gauge, electrical, and thermal sources, have been successfully integrated into long-term condition monitoring protocols \cite{inaudi2010long}.
The richness of spatial-temporal data obtained from these sources lends itself well for use in inversion-based SHM, i.e. given a set or sets of static SHM data, use an inverse methodology in capturing (potential) damage.
This is true in cases where numerical models are available for the problem physics and when they are not (e.g. learned models can act as surrogates when physics-based models are unavailable).

The sheer volume of literature available reporting the successful use of static inversion methods in SHM is formidable. 
However, roughly speaking, static inverse methodologies have been implemented within three areas, (a) point sensing, (b) area sensing, and (c) volumetric sensing.
Holistically, it can be difficult to distinguish between each of these areas; for example, when lower dimensional measurements are extrapolated to characterise damage in area or volume targets \cite{downey2017damage,zhou2012optical}.
One such example is Digital Image Correlation (DIC), where the displacement field at discrete points on a structure is inversely computed via pixel movement and then extrapolated (interpolated) to a full-field, whereby the quality of the computed field is highly dependent on the quality of the contrasting area speckle pattern \cite{lecompte2006quality}.
Similarly, one may consider the distinguishing of static measurement dimensionality as a local-global phenomenon where discrete local changes contribute to the analysis of the global structural system \cite{downey2017damage}.
Lastly, to complicate matters even more, the use of discrete measurements can yield 2- or 3D information -- as in the case of penetrative electrical measurements, where currents diffuse through the entirety of a body \cite{borcea2002electrical}.

Fortunately these, perhaps philosophical realisations, are often washed out via the nature of inversion methodologies themselves.
Pragmatically, at least in the context of SHM, the solution to static inverse problems generally requires a model, either physics based or learned.
As such, the amalgamation or assimilation of data and solution dimensionality is often simply a matter of discretisation or model generation.
In a similar vane, when static inverse problems are ill-posed, solutions generated using lower dimensional data are regularised/biased using prior models consistent with the solution dimension.

Many examples are available in the literature illustrating the efficacy of static inversion methodologies for applications in a suite of SHM implementations.
For example, the use of displacement measurements for capturing SHM causalities in various structural geometries was reported in \cite{zare2016optimization,kefal2016displacement,sanayei2015automated}.
Of note, the specific applications using displacement fields to reconstruction elastic ans elasto-plastic properties (and corresponding damage characteristics) has been the source of significant research \cite{ni2021deep,waeytens2016model,ruybalid2016comparison,karageorghis2016method,mathieu2015estimation,kefal2016quadrilateral,smyl2018coupled,shen2010inverse,xu2010determination,avril2008overview,lecampion2007model,amiot2007identification,hild2006digital,bonnet2005inverse,morassi2004stable,maniatty1989finite}.
In the pervasive case where displacement/strain measurements are discretely measured from strain gauges/fibre optics, inverse methodologies have also been fruitfully employed for damage characterisation, pressure and strain mapping, and shape sensing \cite{gomes2020sensor,an2019recent,colombo2019definition,amanzadeh2018recent,di2015fibre,cerracchio2015novel,gherlone2014inverse,derkevorkian2013strain,weisz2013evaluation,gherlone2012shape,tessler2011real,kiesel2006intrinsic}.
Perhaps one measure illustrating the success of such inverse approaches is highlighted by the recent interest in optimising the related sensing schemes \cite{song2021optimal,esposito2020composite,gomes2020sensor,ostachowicz2019optimization}.

	\label{compbeam}

\begin{figure*}[ht]
    \centering
    \begin{tabular}{cc} 
    \includegraphics[width=0.45\textwidth]{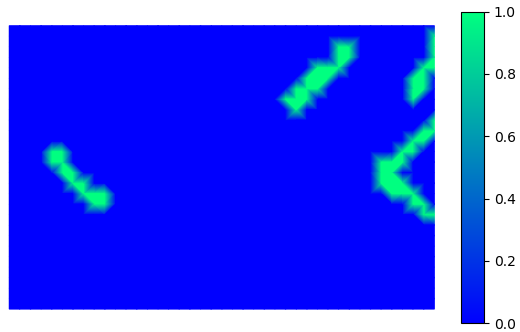} & \includegraphics[width=0.45\textwidth]{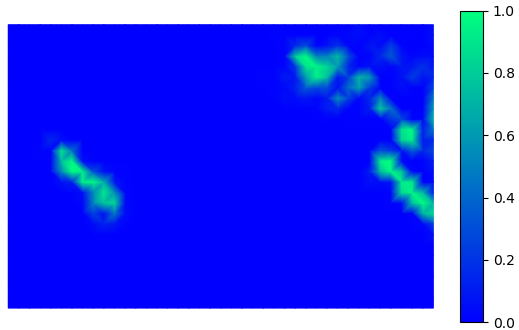}\\
   {a) Simulated shear crack patterns} & {b) Reconstruction of shear crack patterns}\\
    \includegraphics[scale = 0.03]{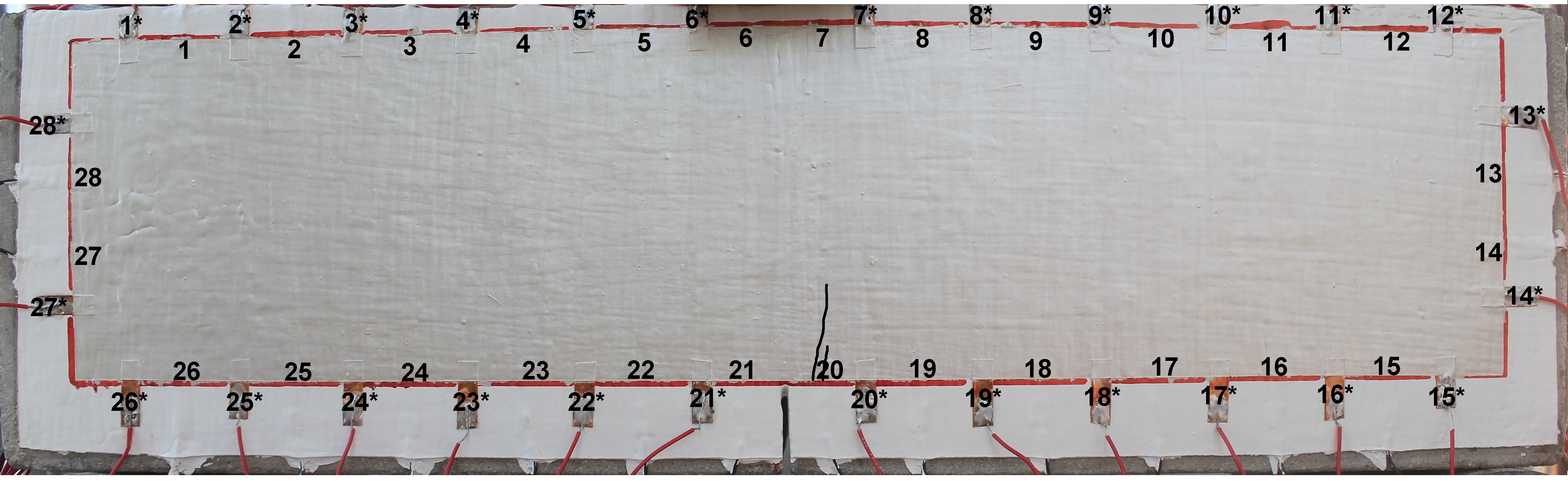} & \includegraphics[scale = 0.45]{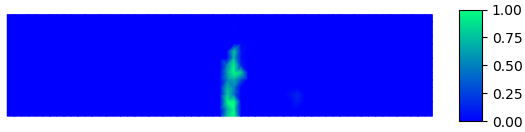}\\
    {c) Experimental image of a flexural crack} & {d) Reconstruction of a flexural crack}\\
    \end{tabular}
        \caption{Reconstructions (right column) report  probabilistic predictions of local flexural and shear cracking in a concrete elements. The colour bars represent the probability of cracks at nodal locations: (a) simulated shear cracking pattern, (b) probabilistic prediction of the shear crack pattern using a convolutional Neural Network, (c) Photo of a flexural crack in an area sensing skin, and (d) probabilistic prediction of the flexural crack using a feedforward neural network.}
    \label{fig:crack reconstruction}
\end{figure*}

In the past decades, the emergence of electrical inverse methods has also proven a viable approach to static condition monitoring \cite{tenreiro2021review,alessandrini2002detecting}.
This family of inversion-based modalities generally utilises three different data sources including capacitative, direct current, and alternating current based measurements.
Capacitative sensing is perhaps the newest of these approaches to sensing, where SHM causalities can be inversely computed using smart bricks \cite{downey2017smart}, area sensors \cite{sadoughi2018reconstruction,kong2018sensing,kong2017large,downey2016reconstruction}, and electrical capacitance tomography \cite{VOSS2019107967,KrzysztofECT}.
On the other hand, owing to its established history in medical and geophysical applications, electrical impedance tomography (EIT, or electrical resistance tomography, ERT) is becoming a well established approach to damage detection, reconstruction, and localisation especially in concrete applications.
For these, EIT presently manifests via two approaches, reconstructing conductivity maps using boundary voltages measured from area sensing skins \cite{smyl2018detection,seppanen2016} or directly imaging a 3D cementitous body \cite{karhunen10ccr}.
Representative 2D EIT reconstructions are provided in Fig. \ref{fig:crack reconstruction} using a machine learned approach (direct reconstruction) for imaging flexural and shear cracks in concrete elements.
To this end, EIT has also been used for characterising area corrosion \cite{seppanen2016} and localising area temperature variations \cite{Rashetnia_2017}.

In summary of this subsection, it is clear that the use of inverse methodologies in static SHM applications is pervasive.
Meanwhile, a number of inversion-based modalities are still emerging -- or are yet to emerge.
Indeed, given the number of potential static data sources available at present, there exist substantial opportunities to investigate or formulate new inversion based modalities.
In the light that some physical models for various underlying physics remain unavailable (either in open source or in general), the use of learned models may bridge this gap.
Lastly, there currently exists tremendous opportunities in the areas of data fusion and joint imaging, which remain predominately unexplored in the area of static inversion based SHM \cite{hauptmann2021fusing}.

\subsection{Dynamical inverse problems in SHM}
The use of dynamical data for monitoring the health and condition of structures is well established \cite{bao2019state}.
For this, a number of data sources are available, for example discrete acceleration, strain, displacement measurements, and recently, coupled electromechanical impedance via piezoelectric transducers \cite{fan2018identification,kim2019electromechanical}.
In the case of typical civil infrastructure \cite{brownjohn2007structural}, the ability to actively excite monitored structures is pragmatically challenging due to the extreme magnitude of the excitation required to attain a distinguishable response.
For this reason, ambient monitoring methodologies have gained significant popularity in recent years \cite{entezami2019structural,ramos2013dynamic,magalhaes2012vibration}.
Irrespective of the dynamical monitoring approach used, extracting dynamical structural properties of interest can be viewed as an ill-posed inverse problem \cite{nagarajaiah2017modeling}.
The ill-posedness of such problems results from a number of actualities, not limited to uncertainties in environmental conditions (wind, temperature, ground conditions, humidity, etc.), traffic, measurement noise, the discrete nature of measurements, material characteristics and numerical modelling error.

One of the most pervasive frameworks used in solving dynamical inverse problems is model updating, which generally aims to match a physics-based model (such as a representative finite element model) to measured dynamical data \cite{worden2008review,friswell2007damage}, commonly using a form of modal analysis \cite{koo2013structural,marcuzzi2010dynamic,pines2006structural}. The physics-based techniques are particularly efficient in providing higher accuracy when testing is restricted.
It is often the case that reconstructing the dynamical SHM properties of interest proves difficult, requiring an innovative approach; some proposed frameworks have included advanced optimisation protocols \cite{hong2010reconstruction} and mode decomposition/superposition \cite{wan2014structural,wan2014structural}.
Alternatively, the use of phase space \cite{paul2017phase}, state space \cite{nichols2003using}, singular value decomposition,
{$\lambda$-curves method \cite{fernandez2017lambda}}\cite{liu2014damage}, auto-regressive \cite{yao2012autoregressive,nair2006time}, Gaussian process \cite{dervilis2016exploring}, and Bayesian/stochastic approaches \cite{ramancha2020bayesian,wan2019bayesian,wan2016stochastic,mao2013statistical} have proven successful.
As noted in the previous subsection, one metric for assessing the progress in this field is the number of works aiming to optimise sensing information, for example in \cite{cantero2020optimal,capellari2018cost,ostachowicz2019optimization,yi2016multiaxial,sun2015optimal,bhuiyan2014sensor,yi2012modified,yi2011optimal,guo2004optimal}.

In the past twenty years, guided wave based modalities have emerged as a viable approach to dynamical inversion based SHM \cite{cantero2021structural,li2019structural,shen2014structural,srivastava2010quantitative,mitra2016guided,eiras2014evaluation,giurgiutiu2005tuned,shin2007improved}.
Common physical manifestations of guided waves include Lamb waves (propagating through thin shell and plate structures) \cite{xu2007single}, Rayleigh waves (surface waves) \cite{wang2004synthetic}, and shear waves.
Generally speaking, SHM systems consist of transducer systems used for actuation and measurement accompanied with an inversion algorithm aiming to reconstruct SHM causalities of interest.
Owing to a number of numerical challenges, conventional solutions to related inverse guided waves problems are generally not feasible \cite{mitra2016guided}.
Consequentially, alternative methodologies have been proposed including, for example, inverse filtering \cite{moll2010time}, reverse time migration \cite{he2020guided,he2019multi}, Bayesian/probabilistic \cite{zhao2020sparse,lu2010evaluation}, amongst emerging inversion approaches.

While the use of dynamical inversion based techniques is well established in conventional monitoring,in some areas (such as guided wave monitoring) it remains in the early stages of development and affords numerous research opportunities. It is worth mentioning that with the advent of modern machine learning methods, we can only anticipate significant advances in forthcoming years as trained networks are now capable of addressing key SHM challenges related to, for example, model error estimation/correction \cite{smyl2021learning} and reducing computational demands associated with many SHM facets \cite{yuan2020machine,finotti2019shm}.

\subsection{Computer vision inverse problems in SHM}

{\color{blue} In addition to the discussed static and dynamic inverse problems in SHM, computer vision-based SHM methods have become an emerging field in inverse engineering problems in approximately the past decade.} Relying on digital images and videos, vision-based SHM techniques enable affordable and rapid structural prognosis. The concept of vision-based inverse problems is straightforward: visual information from the external surfaces of structures is captured through digital cameras, serving as input data for computer vision algorithms in detecting, localising, and quantifying structural damage in a variety of contexts.

Computer vision-based inverse problems can be either static or dynamic in nature. In dynamical environments, a digital camera is treated as a vision sensor measuring dynamic structural responses. Instead of directly capturing the structural vibration through contact-based sensors (e.g. accelerometers), vision-based algorithms can track structural responses through a video stream. For example, in \cite{khuc2017completely}, researchers applied a video feature tracking technique to measure the pixel movements of a steel girder in a football stadium under a service load using a consumer-grade digital camera. These movements were then converted into displacements using a scaling factor.
Similar efforts have been reported in \cite{feng2015nontarget, luo2018robust}. Furthermore, through the usage of cameras as displacement sensors, other key structural features, such as natural frequencies/mode shapes \cite{yoon2016target, yang2017blind}, beam influence lines \cite{dong2019completely}, and bridge cable loads \cite{jana2021computer} have been be estimated. 

In addition to tracking the surface motion, computer vision algorithms can offer rapid and reliable inspections against structural damage such as cracks \cite{liu2019computer}, concrete spalling \cite{gao2018deep}, steel corrosion \cite{shen2018human}, and other structural deterioration \cite{pauly2017deeper}. To make it is viable, researchers develop vision-based algorithms to scan and extract damage-induced visual features either across the entire image scene or within a small predefined image patch (e.g. region of interest) that is prone to structural deterioration. In general, the image-based damage extraction techniques can be categorised as: (a) machine or deep learning-based methods; and (b) non-learning based methods. 

\begin{figure}[h!]
	\centering
	\includegraphics[width=\textwidth]{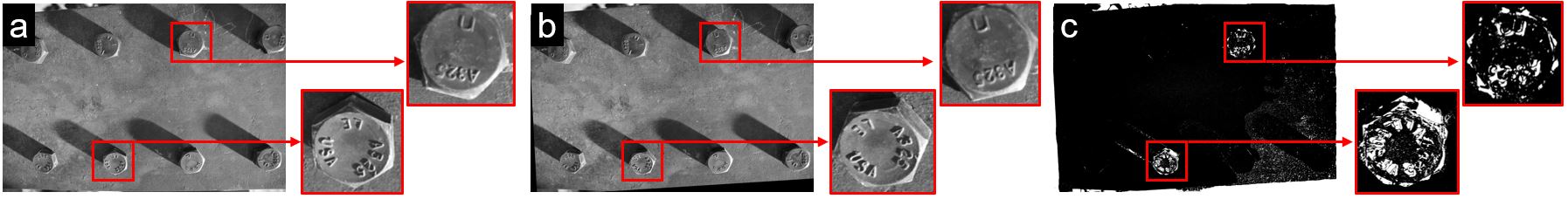}
	\caption{An example of vision-based bolt loosening detection where image (a) and (b) are images of a bolted connection taken at different inspection periods. Two loosened steel bolts are shown in the blow-up figures with counterclockwise rotations in their bolt heads. Using a series of image processing techniques, the differential features caused by the bolt loosening can be identified in c). Detailed discussion can be found in	\cite{kong2018image}.}
	\label{boltloosening}
\end{figure}

The idea of machine or deep learning-based (computer vision) methods is to train a damage detection classifier through an image dataset with pre-labeled structural damage. 
Then, the classifier is applied in characterising structural damage from newly captured images. 
Some of the successful applications include detection of concrete cracks \cite{cha2017deep, dung2019autonomous} and spalling \cite{gao2018deep}, steel cracks \cite{chen2017nb}, bolt loosening \cite{cha2016vision, huynh2019quasi}, steel surface defects \cite{gao2018deep}, pavement cracks \cite{pauly2017deeper}, and complex situations where multiple damage types exist \cite{cha2018autonomous}. In contrast, non-learning based methods can directly pick up image features caused by structural damage, and hence do not require any prior knowledge in training the classifier. For instance, fatigue cracks in steel bridges can be identified through crack breathing behaviour \cite{kong2019non}. Also, loosened bolts in steel connections can be quantified by extracting the differential features provoked by bolt head rotations \cite{kong2018image}. Fig. \ref{boltloosening} illustrates an example by comparing two images of a steel connection at different inspection periods to extract the differential features provoked by the loosened bolts.

Computer vision-based methods can be extended for damage detection and pattern recognition in full scale civil structures. Utilising the platform of unmanned aerial vehicles (UAVs), the on-board UAV camera can rapidly scan the structure including the locations that are challenging to be accessed by traditional contact-based sensors. For example, researchers in \cite{dominici2017uav, choi2018computer} applied UAVs and vision algorithms to leverage effective approaches for post-earthquake building safety inspections. Similar efforts have been reported for inspecting dams \cite{buffi2017survey, khaloo2018utilizing}, tunnels \cite{ozaslan2015inspection}, and railways \cite{banic2019intelligent}. Other researchers adopted satellite images to examine damage status over a larger scope of work (i.e., multiple buildings at the community level) after natural disasters such as flood, earthquake, volcanic eruption, hurricane, and wildfire \cite{shao2020bdd, gupta2021rescuenet}.

UAV platforms are also capable of collecting a large volume of images of civil structures under different camera angles through automated route planning. Such an advantage can be further enhanced by a computer vision workflow, termed photogrammetry, for the purpose of reconstructing the 3D model of the structure. Relying on structure-from-motion with multi-view stereo (SfM-MVS) algorithms \cite{furukawa2015multi}, photogrammetry technique can create a 3D structure model (e.g. a 3D point cloud) based on 2D images. Photogrammetry leverages several potentials in inverse SHM problems. 
For instance, in recent work, \cite{braun2019combining} created a dense point cloud of a building in the construction site based on UAV images. Then the point cloud was integrated with a building information modelling (BIM) model for labelling the structural components in original UAV images. In \cite{stavroulaki2016modelling}, the researchers utilised the dense point cloud to assist the creation of a finite element model of a masonry bridge. The authors argued the benefit of using the point cloud was twofold: the point cloud depicted the accurate geometric information of the bridge and offered the results of bridge crack distribution. In \cite{zollini2020uav}, researchers leveraged the photogrammetry workflow to find concrete cracks and spalling of a highway bridge. Lastly, in the context of the 2021 Hernando de Soto Bridge incident \cite{leishear2021bridge}, where a large crack was discovered in a "fracture critical" element by a private engineering firm -- yet, previously identified approximately two years earlier by a local operating a commercially available drone: the use of coupled UAV/computer vision approaches to SHM may be more valuable than ever.

\subsection{Digital twins and outlook}
As this section has illustrated, the use of inverse methodologies in SHM is both well established and an area of active development.
With the rapid digital transformation of structural assessment and infrastructure asset management, the emergence of numerous digitally-inspired technologies will play a key role in the future trajectory of inverse problems in SHM.
At the forefront, \textit{digital twins} have been the source of increasing research and industrial interest in recent years \cite{wagg2020digital}.
While the scope of digital twins' applications spans beyond SHM alone, its basic aim is to provide information on the current or future state of an asset by combining real-time data and a physical/data-driven model offers many potential avenues for engagement with the inverse problems community.
Nonetheless, in specifically considering a classical SHM application, such as damage localisation \cite{seshadri2017structural}:
developments stemming from the inverse community including, for example, state estimation \cite{evensen2018analysis,iglesias2013ensemble,seppanen2001state}, uncertainty/model error approximation/compensation \cite{lunz2021learned,smyl2021learning,huttunen2007approximation}, regularisation \cite{benning2018modern}, and model reduction \cite{lieberman2010parameter}, have excellent potential for enriching or enhancing digital twin frameworks.

As a whole, the future outlook for the integration and advancement of inverse methodologies in SHM is very bright.
Indeed in the past 20 years, we have seen an exponential increase in high-performance computing and graphical processing units development and assimilation into modern civil and mechanical engineering applications \cite{nagel2019high,gaska2018high}.
Coupled with powerful inverse frameworks for large-scale problems (e.g. Krylov solvers \cite{saibaba2013flexible} and distributed computing \cite{yang2017coded}) and machine learning \cite{lucas2018using}, we can only expect a steady increase in (a) the breadth of inverse problems the SHM community is able to address and (b) and an evolution in innovative inversion-based approaches to solving increasingly challenging SHM problems.

\section{Smart materials and structures}\label{smart_materials}

\emph{Smart}, \emph{self-sensing} materials have received immense attention in recent decades \cite{aviles2018piezoresistivity, han2015intrinsic, han2020review}. A material is said to be self-sensing if it exhibits a property change in response to external stimuli. These materials are able to intrinsically report on their health or condition in a spatially continuous manner and with less hardware/instrumentation burden than traditional point-based sensing technologies (e.g. strain gauges, piezoelectric patches, accelerometers, etc.). In structural engineering, external stimuli are often mechanical effects such as deformation, damage, or loads. Hence, integrating smart materials into next-generation structures may allow for unprecedented health monitoring and diagnostics. {A discussion on smart materials in structural contexts may therefore also be considered as a subset of the preceding discussion on SHM (see \S\ref{shm}). Nonetheless, we will treat it as a distinct topic in the forthcoming due to its unique inverse problems.}

Although self-sensing is an umbrella term encompassing many different physical effects, the \emph{piezoresistive effect} has perhaps received the most attention to date (see recent reviews \cite{aviles2018piezoresistivity, forintos2019multifunctional, tallman2020structural}). Piezoresistive materials are so-named because they exhibit a change in electrical conductivity (or its inverse, resistivity) upon deformation. 
This means that every point of the material is capable of relaying information on its mechanical state. 
Damages such as voids, ruptures, or fractures can also be detected since the removal of material represents a conductivity loss. 
But spatially-continuous piezoresistive sensing presents two challenges of relevance to inverse problems: (a) It is not practical to instrument electrodes at every point on a structure, which means that it is necessary to deduce conductivity distributions from a finite set of measurements. And (b), even if we could recover a spatially continuous mapping of the conductivity, electrical properties are of little consequence to the structural engineer. We would much rather know the underlying mechanics that give rise to an observed conductivity distribution. We will address both of these inverse problems. First however, a brief summary of the physics and piezoresistive modelling approaches is given. It will be seen later that modelling techniques are essential to solving inverse problem (b) above.

\subsection{Piezoresistive nanocomposites}\label{piezoresistive_nanocomposites}

Many materials intrinsically exhibit piezoresistive properties. For example, carbon fibre-reinforced polymers (CFRPs) are well known to change conductivity when loaded elastically \cite{wang2013through,koecher2015piezoresistive}. Here however, we will instead focus on materials that have been \emph{engineered} to be self-sensing; that is, an additional constituent has been added to the material system without which it does not exhibit piezoresistivity. This is most commonly done by adding a conductive phase to a non-conductive matrix such as polymers (including structural polymers like  epoxy vinyl ester \cite{ku2018selective}, polymeric thin films for use as sensing skins \cite{hou2007spatial}, laser-induced graphene inter-layers in continuous fibre composites \cite{groo2021fatigue, groo2021damage}, and even polymer binders in energetic materials \cite{talamadupula2021piezoresistive, talamadupula2020mesoscale}), cements \cite{yoo2018electrical}, or ceramics \cite{ricohermoso2021piezoresistive}. Electrical transport is then a consequence of percolation -- the composite conducts electricity when enough fillers have been added to form an electrically connected network. Because the percolation threshold decreases with aspect ratio, fillers with ultra-high aspect ratios like carbon nanotubes (CNTs) are popular. There are considerable challenges associated with manufacturing CNT-based nanocomposites such as achieving good dispersion -- ultra-small fillers such as CNTs have a tendency to agglomerate which can degrade the mechanical properties of the composite. But manufacturing is outside of the scope of this manuscript and is well-covered elsewhere \cite{mccrary2012review,ma2010dispersion}.

Considerable effort has also been devoted to the development of piezoresistivity models -- computational and/or analytical means of predicting how conductivity changes for a prescribed strain. These  efforts can be broadly categorised as (a) equivalent resistor network models \cite{gong2014mechanism,hu2008tunneling,lee20152d}, (b) computational micro-mechanics models \cite{talamadupula2021piezoresistive, alian2019multiscale, ren2015modeling, talamadupula2021statistical}, or (c) homogenised macroscale models \cite{cattin2014piezoresistance,tallman2013arbitrary,koo2020higher}. In (a) equivalent resistor network models, high aspect ratio fillers like CNTs are represented as either straight or wavy/curved sticks in a micro-domain \cite{gong2014mechanism,hu2008tunneling,lee20152d}. These sticks are discretised into resistors based on the length, diameter, and conductivity of the fillers, and junctions between nearby fillers are discretised into resistor elements based on the equivalent resistance felt by an inter-filler tunnelling electron \cite{bao2012tunneling}. The conductivity of the nanocomposite can then be calculated from the overall resistance of the discretised nanofiller network and the dimensions of the micro-domain. For a given deformation of the micro-domain, the translation and rotation of the fillers can be calculated by treating them as rigid inclusions \cite{taya1998piezoresistivity}. Post-deformation, the conductivity of the micro-domain is recalculated thereby allowing piezoresistive properties to be predicted.

Second, (b) computational micro-mechanics models use computational means to simulate both phases of the composite -- the non-conductive matrix and the conductive fillers \cite{talamadupula2021piezoresistive, alian2019multiscale, ren2015modeling, talamadupula2021statistical}. Because of this, and unlike equivalent resistor network models, computational micro-mechanics models can incorporate nuanced effects such as nanofiller-to-matrix debonding, nanofiller deformation, etc. A limitation of this approach is computational expense due to individual nanofillers and the enveloping matrix being explicitly simulated. It is therefore difficult to scale these models to structural levels. 

And third, (c) homogenised macroscale models describe the conductivity of the nanocomposite as an analytical function of the strain state without simulating individual fillers. Rather conductivity/resistivity-strain are coupled through analytical functions based on excluded volume theory \cite{tallman2013arbitrary} or via `constitutive' tensors (note that these are not proper constitutive relations because conductivity and strain are not energy complements) \cite{koo2020higher,zhao2018characterization,gruener2017characterization}. This approach is therefore less computationally expensive than equivalent resistor network and computational micro-mechanics models. Importantly, homogenised approaches can be readily integrated with structural analysis tools such as the finite element method, thereby allowing for macroscale piezoresistive analyses. {As will be discussed in \S\ref{smart_materials}\ref{piezoresistive_inversion}, this allows for strain recovery via piezoresistive inversion.} Despite these advantages, analytical models suffer from having to make assumptions regarding average inter-filler spacing, average orientation of nanofillers, and the need for calibration data. 

\subsection{Conductivity imaging via EIT/ERT}

{As discussed in \S\ref{shm}\ref{static_inverse_problems_in_shm},} conductivity imaging modalities such as EIT (or DC resistivity imaging via ERT) have been explored for health monitoring in structural materials \cite{tallman2020structural}. EIT is a natural complement for piezoresistive materials because it allows for the spatial localisation of damage and the \emph{mapping of deformation and strain.} There are several factors that make the EIT interesting to pair with self-sensing materials: (a) This combination allows for sub-surface strain imaging. That is, a myriad of techniques exist for monitoring surface strains such as strain gauges, DIC, holographic methods, etc. Tools for studying sub-surface strains, however, are much more limited, often require ionising radiation, and can be costly (e.g. volumetric strains via X-ray digital volume correlation \cite{bay2008methods}). (b) Traditional `limitations' of EIT can be strengths in structural applications. For example, the simplest implementation of EIT favours spatially smooth solutions. This is obviously undesirable when imaging discontinuous features with distinct boundaries such as organs in a living body. Strain fields, however, are often smoothly varying thereby playing into EIT's strengths. And (c), because piezoresistive materials are engineered, we can leverage our knowledge of their piezoresistive response to build bounds into the EIT inverse problem. Two examples of strain imaging via EIT in self-sensing polymeric composites are shown in Fig. \ref{smart_materials_EIT_piezo_inversion_fig}. 
Both examples leverage knowledge such knowledge of conductivity change bounds to build constraints into the solution space.
{Beyond these advantages for structural imaging, the pairing of EIT/ERT with self-sensing materials may also have keen, as-of-yet unrealized potential for extreme loading (particularly self-sensing energetic materials \cite{talamadupula2021piezoresistive, talamadupula2020mesoscale}). That is, it may be illuminating to image energetic materials as they detonate. However, most prevailing imaging modalities have not the temporal resolution to capture these fleeting moments and require hardware that is too expensive to risk damaging. EIT, on the other hand, can have ultra-high temporal resolution (on the order of hundreds of microseconds for optimized systems \cite{dowrick2018phase}) and uses only low-cost hardware (i.e. sacrificing the hardware during a detonation is of virtually no financial consequence). Thus, there may be much exciting potential for overlap between smart materials + EIT and extreme loading as described in \S\ref{extreme_loads}.}

There have been many studies on the topic of EIT and piezoresistive materials. A few representative examples are summarised in this manuscript, but interested readers are directed to a recent review for a more in-depth discussion \cite{tallman2020structural}. Some of the first work in this area made use of self-sensing nanocomposite \emph{sensing skins} produced by a layer-by-layer fabrication technique \cite{hou2007spatial, loh2009carbon}. These skins were applied to substrates, and EIT was then used to identify and spatially localise deleterious effects including mechanical etching, pH exposure, impact damage, and residual strains from impacts. Further on the topic of polymer-based composites, thin pressure sensors were produced by embedding a non-woven textile modified with CNTs into a soft elastomer. EIT was then employed to visualise various pressure distributions including non-uniform distributions \cite{dai2019large}. EIT was also recently applied for damage detection in a ceramic-based composite that was modified with micron-sized waste-iron particles \cite{nayak2019spatial}. And as a final representative example, self-sensing cement composites have been produced by spray-depositing a CNT-modified latex on the aggregate phase. EIT was then successfully used to localise various damages in the material \cite{gupta2017self}.

\subsection{Piezoresistive inversion}\label{piezoresistive_inversion}

Even though solving the EIT inverse problem allows for the spatially continuous visualisation of mechanical effects in these materials, this poses an obvious problem -- structural engineers are generally not interested in conductivity. They would much rather know the spatially-varying components of the strain tensor, stress tensor, or precise damage characteristics since these factors drive structural analyses and health assessments. Recalling also that various macroscale piezoresistivity models exist as described in \S\ref{smart_materials}\ref{piezoresistive_nanocomposites}, we can formulate another inverse problem as follows: ``for a given EIT conductivity distribution (or, more directly, for a given set of EIT boundary data) and with an accurate model of conductivity-strain or conductivity-damage coupling for a particular material, can the precise mechanics of the system be recovered?''

\begin{figure}[ht]
	\centering
	\includegraphics[width=0.85\textwidth]{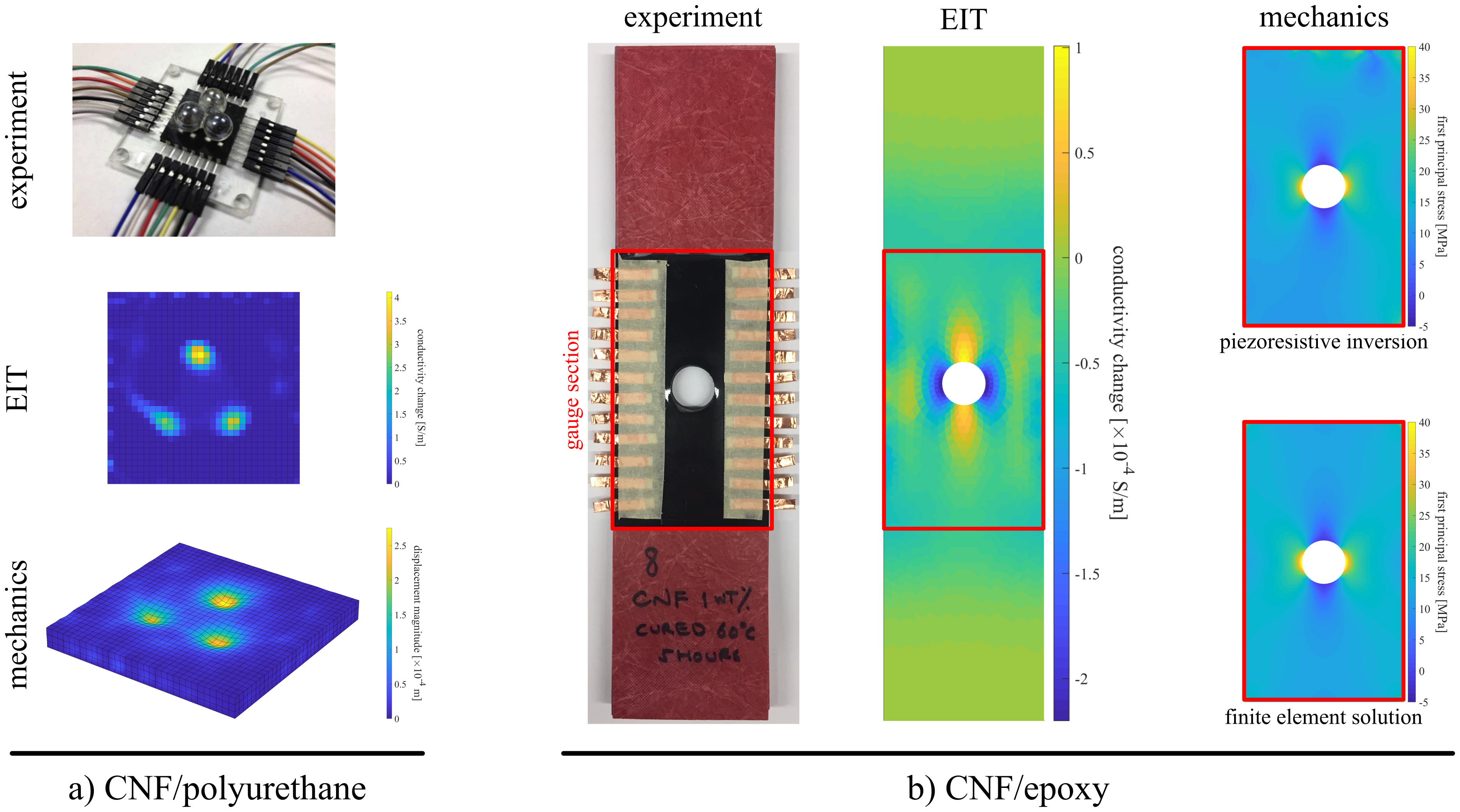}
	\caption{Examples of EIT and piezoresistive inversion applied to self-sensing nanocomposites. a) A soft CNF/PU is deformed by rigid, non-conductive indentors \cite{tallman2017inverse}. EIT is then used to image the deformation-induced conductivity changes, and piezoresistive inversion is used to recover the displacement field (multiplied by a factor of five for ease of visibility). b) A hard CNF/epoxy is loaded in tension with a stress raiser at its centre \cite{hassan2020failure}. EIT is again used to image the conductivity change. Lastly, piezoresistive inversion is used to recover the underlying displacement field. With knowledge of the material's elastic properties, strains and stresses can be spatially mapped. The first principal stress of the guage section is shown here along with comparison to a traditional FEM solution for validation.}
	\label{smart_materials_EIT_piezo_inversion_fig}
\end{figure}

Both of these piezoresistive inversion problems\textemdash strain recovery and damage recovery\textemdash are challenging. The former is challenging because, under ideal circumstances, we seek six components (in 3D) of a strain tensor from a scalar conductivity field. Prospects can be improved somewhat by instead seeking the displacement field (i.e. three unknowns in 3D) from the conductivity data and making use of reasonable assumptions (e.g. plane strain and plane stress) where applicable, but the displacements-from-conductivity inverse problem is nonetheless under-determined. The challenge is exacerbated by the fact that circumstances are never truly ideal; conductivity and conductivity changes are not exactly isotropic, and EIT cannot image individual components of a conductivity tensor. Even for the case of damage imaging, an accurate model of material breakage-induced conductivity change is needed. For simple nanofiller/matrix phase nanocomposites (i.e. without reinforcing fibre), material breakage can be treated as a cessation of conductivity. For more complicated material systems such as nanofiller-modified continuous fibre composites, however, material breakage-induced conductivity changes must account for factors such as anisotropy and residual post-damage conductivity due to inter-laminae contact. And even if a suitable damage model is developed, the inverse solver needs to be capable of reproducing potentially complex damage shapes that are not readily amenable to parameterisation. Both strain and damage recovery are additionally hampered by the fact that EIT does not produce accurate conductivity distributions in an absolute sense or a spatial sense.

Despite these challenges, the piezoresistive inversion problem has received some attention. An initial effort used an analytical inversion framework predicated on iteratively minimising the $l_2$-norm of an error vector between a predicted and observed conductivity distribution \cite{tallman2016inverse}. Although this work was entirely computational and limited to simple deformations and infinitesimal strains, it nonetheless demonstrated that piezoresistive inversion was possible. The next work in this area used EIT to image strain-induced conductivity changes in a carbon nanofiber (CNF)-modified polyurethane (PU) composite \cite{tallman2017inverse}. Three marbles (i.e. comparatively rigid non-conductive indentors) were pushed into the CNF/PU as EIT measurements were taken. A similar analtyical approach was used to reproduce the displacement field. 
Two important factors differentiated this study -- experimental validation of piezoresistive inversion and successful application to materials undergoing finite strains. 
Later works looked at utilising metaheuristic algorithms for solving the strains-from-conductivity inverse problem in a CNF-modified epoxy \cite{hassan2020failure,hassan2020comparison}. The CNF/epoxy was moulded in the shape of a plate with a hole and loaded in tension, causing strain concentrations in the vicinity of the hole. Genetic algorithms, simulated annealing, and particle swarm optimisation were explored because it was observed that the analytical formulation failed to converge to the physically correct solution for this more complex loading state. It was found that genetic algorithms performed the best for this inverse problem, but all methods compared favourably to DIC experimental validation. Because epoxy is relatively brittle, these studies were necessarily limited to infinitesimal strains. Despite the successes of the preceding studies, they were all limited to electro-mechanically isotropic materials. Translating these capabilities to electro-mechanically anisotropic materials remains a daunting challenge. Fig. \ref{smart_materials_EIT_piezo_inversion_fig} summarises results for strain recovery via piezoresistive inversion.

Some work has also been done for damage recovery via piezoresistive inversion. Recall that a key challenge with the damage recovery inverse problem is shape parameterisation. To that end, these preliminary studies have considered relatively simple damage cases. For example, various machine learning algorithms were used to categorise three damage conditions in a self-sensing bone cement directly from EIT boundary voltage data \cite{ghaednia2020machine}. {Thus, self-sensing materials also have overlap with machine learning methods (see \S\ref{MLInversion}).} Bone cement is poly(methyl methacrylate) (PMMA) used to facilitate robust contact between an othopaedic implant (e.g. a total joint replacement) and hard bone. Failure of the PMMA interface is often difficult to detect via radiographic imaging, hence the motivation for alternative diagnostic tools. This clarification aside, the parameterisation was relatively simple in this case -- only four distinct states were possible (three damaged plus one healthy). The combination of EIT, self-sensing PMMA bone cement, and machine learning allowed for correct damage classification with over 90\% accuracy. In another study, image recognition-based machine learning was used to identify, size, and localise through-hole damage to a self-sensing composite plate \cite{zhao2020spatial}. The image recognition algorithm was trained using computationally generated EIT images on a simple square domain punctured by a random number of randomly sized circular holes. The trained network was able to adeptly predict through-hole size, location, and number from EIT conductivity maps with good accuracy -- likely better accuracy than human interpretation of EIT images. But this again utilised a very simple damage parameterisation (i.e. needing only to predict hole number, radius, and in-plane coordinates). As a final example, a recent study looked at delamination shaping from EIT images in CNF-modified glass fibre/epoxy laminates \cite{hassan2020damage}. In this study, delaminations induced by low-velocity impacts were parameterised as ellipses of unknown major and minor axes and centred at unknown in-plane coordinates. A genetic algorithm was used to inversely determine these parameters by minimising the $l_1$-norm of the difference between experimentally collected EIT boundary voltage data and boundary voltage data predicted by a computational model of the damaged domain. Destructive analyses of the post-impacted laminates revealed that the GA-predicted damage state much more closely matched the actual delamination size and shape than the EIT conductivity images. This third example of damage recovery is particularly noteworthy because it represents a much more realistic damage state.

In summary, this section has looked at smart, self-sensing materials from the perspective of structural inverse problems. Two noteworthy inverse problems were discussed -- the EIT inverse problem and the strain/damage recovery problem. The former has been extensively researched in other fields (e.g. \cite{holder2004electrical}). The latter, however, is much more recent and has only been the subject of a few precursory studies. Much work remains to be done regarding the inversion of electrical data to obtain underlying mechanical effects. {Nonetheless, it can be seen that the field of smart materials in structural applications has much potential and cross-cutting overlap with other topics of this article including extreme loading, SHM, and machine learning.}

\section{A look forward: machine learning and education}\label{MLInversion}

This paper has evidenced the pervasiveness of inverse problems and methodologies used in the field of structural engineering. Yet, the following remains to question: "what is guiding the future trajectory of inversion in structural engineering?" 
In this section, we examine what we foresee as the two most influential future areas in addressing this question: machine learning and education.

\subsection{Machine learned inversion}
Many areas of structural engineering rely heavily on applied mathematics and science -- from the integration of material models within finite element frameworks to experimental measurement of structural response excitation.
From a broad scientific perspective, there are two paradigms to research, either (a) the Keplerian paradigm (data driven, obtaining discoveries via data analysis) or (b) the Newtonian paradigm (first principles, discovery through fundamental principles) \cite{weinan2021dawning}.
Without question, structural engineering research uses both principles.

Often, first principles approaches are manifested via partial differential equations (PDEs) and their analytical or numerical solutions.
However, both in research and practice, obtaining solutions to PDEs can be infeasible or intractable, for example owing to computational demand, a dearth in available numerical regimes and/or the ``curse of dimensionality'' \cite{han2018solving}.
In such cases where engineering problems are governed by such PDEs, solving a related inverse problem using a conventional methodology would also be a dubious task.
When faced with this situation, we are most likely constrained to adopting a Keplerian approach.

Unmistakably, machine learning has provided the science and engineering communities with a powerful tool for data-driven analysis, prediction, assessment, and significantly more.
Structural engineering research has also greatly benefited from machine learning, especially in the areas of performance assessment \cite{sun2020machine}, SHM \cite{farrar2012structural} and analysis of various structural phenomena \cite{lee2018background}.
Yet, significantly less attention has been paid to the use of machine learning for solving inverse problems in structural engineering outside of areas such as SHM and NDE (as identified in \S\ref{shm}).

Exemplifying this reality, structural design highly under utilises  machine learning and data science. Design is traditionally associated with an iterative nature, in which various structural concepts are tested conceptually until a prevailing option is identified which adheres to constraints initially identified. The nature of this iterative design process has been understood in the past using both positivist, pragmatic and post-modernist epistemologies ranging from Simon’s ``science of the artificial'' \cite{simon1996sciences}, Sch{\"o}n’s ``reflexive practice'' \cite{schon1983reflective} and Buchanan’s ``placements for contextualisation'' \cite{buchanan1992wicked} respectively. The complexities involved in design from satisfying conflicting demands to exercising appropriate judgements is hence often attributed as being an innate human skill \cite{Addis1990}.

However, recent developments successfully solving inverse problems using data-driven approaches suggest that such methodologies could also be incorporated within structural design \cite{arridge2019solving}. This challenges the notion that design is an exclusive human ability, a development which mimics the success of self-driving cars through data-driven approaches \cite{badue2020self}. What is typically considered ``intuition'' or ``engineering judgement'' may in fact be recalling ``data-points'', committed to long-term memory through the process of experience, that are suitable solutions based on the unique set of constraints one is presented with \cite{simon1996sciences}. {Indeed, supervised machine learning was originally understood to create mapping functions that correlate a set of inputs (a specific set of constraints for a design situation) to associated outputs (viable structural solutions) by learning from a given data set (experience) \cite{mitchell1997machine}. We note that this view has significantly evolved with the onset of deep-learning, as well as the use of unsupervised and reinforcement learning regimes.}

We therefore believe that this inverse problem perspective evidences the need for research on the development of machine learning tools for structural design. Data-driven approaches are often dismissed because they are ''black-boxes'' which lack a scientific rationale for their outputs. However for design, this criticism might be unwarranted, and instead highlights the dogmatic concentration within academia on exclusively solving forward problems, which due to their well-posed nature, lend themselves to engineering science thinking. Recently, researchers have started investigating data-driven models. In one example, neural networks and clustering were applied to form a bridge and navigate the design space and shortlist viable and fitting solutions \cite{karla,ZHENG2020103346}. Other examples include researchers applying machine-learning models to build suitable structural predictors for conceptual design related to building massing \cite{wang2018conceptual}.

The former developments are, in our view, only a starting point for the employment of machine learned models used in solving inverse problems in structural engineering as a new horizon of data-driven approaches emerge.
Especially, for the cases involving intractable forward problems, model reduction techniques have been promising \cite{lipponen2018correction,kaipio2006,Kaipio2007,arridge2006a}, but these are either difficult to design by hand or are restricted by overly simplistic assumptions. Here, data driven approaches are a powerful alternative to compensate for modelling errors \cite{smyl2021learning,lunz2021learned,koponen2021model} or reducing computational cost of iterative optimisation schemes by model approximations \cite{smyl2021efficient,banert2020data}. Finally, we note that recent developments in geometric learning extend deep networks on Euclidean meshes to general meshes, such as finite elements, by a embedding them into graph structures essentially utilising the underlying geometry \cite{wang2019dynamic,herzberg2021graph}. This opens the possibility to extend many data driven approaches to complex structural problems.

\subsection{Inverse methodology in structural engineering education}
Over the last 100 years, engineering education has experienced a number of fundamental shifts: a shift in focus away from design to engineering science in the 1960s, the rise of outcome-based accreditation in the UK and USA in the 1990s, along with a re-emphasis of teaching design through capstone projects in the 2000s \cite{froyd2012five}. There also is the continued tension between teaching graduates both the technical knowledge as well as the interpersonal skills demanded from industry deemed necessary to become effective designers \cite{crawley2007rethinking}. To this day, there exists the debate on how to find the required balance between \textit{knowledge-that} and \textit{knowledge-how} as identified in \S\ref{invdes} \cite{bulleit2012makes}. 
For civil and structural engineering disciplines, we believe that some of these challenges might be addressed by communicating the existence of inverse problems, their pervasive occurrences as shown by the sections above, and teaching the various methods and techniques for solving such problems.

The dominance of engineering science within university curricula today, which primarily focuses on identifying and solving forward problems, might unintentionally generate the wrong supposition that all problems in engineering are well-posed with idealised assumptions. Without adequately addressing the existence of inverse-problems, and their distinctively qualitative differences with forward problems, it is easy to mistakenly assume that, for example, structural design is the application of such ``forward problems''. However, the idea that engineering is simply ``putting theory intro practice'' \cite{Addis1990} or ``applied science'' \cite{koen2013debunking} has been strongly argued against by numerous engineers and philosophers \cite{Koen2003, goldman2004we, dym2005engineering, doridot2008towards}. The challenge which students face when dealing with real-world design problems might be accounted for by the fact that during the majority of their engineering education, they might lack a conceptual framework to adequately demarcate design from analysis. Similarly, in order to engage with other promising fields of structural engineering, such as structural health monitoring, blast engineering and smart materials, educating students on inverse problems is crucial.

Hence, a possible improvement for current civil and structural engineering curricula is introducing students to the existence of forward and inverse problems, how they relate to one another and provide examples where each type of problem arises and how to solve them. As identified previously, this will also potentially require teaching students a host of new skills, especially if data-driven models continue to be effective tools for solving inverse problems, as is the case in structural health monitoring and increasingly likely in design. More importantly, especially when taking into consideration the recent developments in data science ranging from CNNs \cite{rawat2017deep}, transformers \cite{vaswani2017attention} and graph neural networks \cite{wu2020comprehensive}, there exists a vast spectrum of knowledge and applications we may not even be aware of. 

As a matter of fact, in terms of research, we are perhaps faced with a unique situation in academia today. Although only time will tell, it could be argued that similar to how the ``invention'' (or discovery) of calculus in the 18th century was instrumental in providing us necessary tools for solving forward problems, resulting in material models and PDEs which allow the creation of complex finite-element methods, so too might the rise of machine-learning and data-science, which is only now starting to gain serious attention in mathematics \cite{weinan2021dawning}, allow a more rigorous treatment of solving inverse problems. By realising the pervasiveness of inverse problems in structural engineering, but also the fundamental differences with forward problems, there is potentially a vast, untouched and exciting realm of research which awaits.

\section{Conclusions}
This article aimed to demonstrate that numerous structural engineering sub-fields may be either fundamentally or partially viewed as inverse problems.
It was shown that this concept is well accepted in, for example structural health monitoring; however, sub-fields such as structural design are not commonly (formally) defined as inverse problems.
We argue that, by shifting this paradigm in structural engineering academia and industry, we may collectively capitalise from the rich methodologies and approaches already established in the inverse problems community.
This beneficial relationship between structural and inverse communities is expected to pay exponential dividends as new tools, such as machine learned models, emerge and develop -- offering new opportunities for solving previously inaccessible, intractable, and/or unforeseen structural challenges.

\section*{Acknowledgements}
AH is funded by Academy of Finland Project 338408 and Project 336796 (Finnish Centre of Excellence in Inverse Modelling and Imaging, 2018--2025).
DS was supported by Engineering and Physical Sciences Research Council Project EP/V007025/1.
DL was supported by the National Natural Science Foundation of China under Grant 61871356.

\bibliography{bibliography}

\end{document}